\documentclass[12pt]{article}

\usepackage{amsmath}
  \usepackage{paralist}
  \usepackage{graphicx}
  \usepackage{amssymb}
\usepackage{xcolor}

\numberwithin{equation}{section}

\def \beq{\begin{equation}}
\def \eeq{\end{equation}}

\begin{document}
\title{Weak asymptotic solution of one dimensional zero pressure dynamics system in the quarter plane}

\author{Kayyunnapara Divya Joseph$^1$ \\
}

\maketitle

\begin{abstract}
In this paper we study a system of equations which  appear in the
modelling of many physical phenomena. Initially this system appeared in description of the large scale structure formation.
Recently it is derived as a second order queueing model. We construct weakly asymptotic solutions of the initial boundary value problem for the system and interaction of waves  
 in the quarter plane $\{(x,t): x>0,t>0\}$ with boundary Riemann solution centered at $(0,0)$ and Riemann solution centered at a point $(x_0,0), x_0>0$. \\
{\bf  AMS Subject Classification:} {35A20, 35L50,35R05} \\
{\bf Key Words and Phrases:}{ Weak asymptotic solution, Adhesion approximation,  Inviscid, Hyperbolic} \\
{ \bf E-Mail: } divediv3@gmail.com
\end{abstract}

\numberwithin{equation}{section}
\numberwithin{equation}{section}
\newtheorem{theorem}{Theorem}[section]
\newtheorem{remark}[theorem]{Remark}

\section{Introduction}
The following non strictly hyperbolic system of conservation laws,
\begin{equation}
\begin{aligned}
u_t+(u^2/2)_x=0,\,\,\,
\rho_t +(\rho u)_x =0,
\end{aligned}
\label{1.1}
\end{equation}
appear in many applications. For example when one studies the plane
wave solution or radial solution of multi-dimensional the zero pressure gas dynamics system
\[
u_t +(u.\nabla)u=0,\,\,\,\rho_t +\nabla(\rho u)=0
\]
we arrive at \eqref{1.1}. Here $u=(u_1,...u_n)$ is the velocity, $\rho$ is the density, $x\in R^n$ space variable, $t>0$
is time, and $\nabla$ is gradient with respect to space. The space dimensions $n=1,2,3$ are physically important and was first derived to describe formation and evolution  of large scale structures such as galaxies. The velocity $u$ and the density $\rho$ 
describe the evolution of the particles. These particles move with constant
velocities until they collide and at collision, the colliding particles stick
together and form a new huge particle, $u$ is bounded pointwise
and $\rho$ is a measure \cite {g1,g2,GSS}.

There are more recent applications of \eqref{1.1} involving kinetic models of stochastic production flows such as the flows of products through a factory or supply chain. These kinetic models when expanded into deterministic moment equations leads to the system \eqref{1.1} see \cite{a1,a2,FGH15}.

The first equation in \eqref{1.1} says that
$u$ is constant along the characteristics. When the 
characteristics cross, the solution is multi-valued and shocks develop, the solution  is understood in a weak sense. Weak solutions are not unique and in order to select the physical solution, one method is to use the adhesion approximation,
\begin{equation} 
\begin{aligned}
u_t+(u^2/2)_x=\frac{\epsilon}{2} u_{xx},\,\,\,
\rho_t +(\rho u)_x =0.
\end{aligned}
\label{1.2}
\end{equation}
or the modified adhesion model,
\begin{equation} 
\begin{aligned}
u_t+(u^2/2)_x=\frac{\epsilon}{2} u_{xx},\,\,\,
\rho_t +(\rho u)_x =\frac{\epsilon}{2}\rho_{xx},
\end{aligned}
\label{1.3}
\end{equation}
see \cite{g1,g2, GSS}. 
The significance of the adhesion approximation \eqref{1.2}  is to select  the right physical solution for the inviscid system \eqref{1.1}, in the vanishing viscosity limit ie, as $\epsilon$ goes to zero. 

Many authors have worked on the initial value problem for the inviscid system  \eqref{1.1} using
different methods. LeFloch \cite{le1} constructed the solution using the Lax formula  and  justified the terms having form $ \rho u$  using the Volpert product \cite{v1} in the equation.  He also proved that for the Riemann
problem, there exist several solutions. The solution constructed in  \cite{j1,j2}
using the modified adhesion approximation  \eqref{1.3} gives the correct physical solution which coincides with the solution constructed in  \cite{le1}  using the Lax formula.

Here we consider the initial boundary value problem, and construct solution to the inviscid system \eqref{1.1} in $x>0,t>0,$  with initial conditions 
\begin{equation}
u(x,0)=u_0(x),\,\,\,\rho(x,0)=\rho_0(x),\,\,\,x>0
\label{1.4}
\end{equation}
with Dirichlet boundary conditions at $x=0$,
\begin{equation}
u(0,t)=u_B(t),\,\,\,\rho(0,t)=\rho_B(t),\,\,\, t>0.
\label{1.5}
\end{equation}
Existence and uniqueness of the solution to \eqref{1.1} with strong form of initial and boundary conditions \eqref{1.4}, \eqref{1.5} is not  generally possible.   In comparison to the initial value problem treated in \cite{ j1,  le1} the boundary value problem is less understood. We cannot prescribe the boundary conditions \eqref{1.5} for $u$ and $\rho$  in the strong sense
because the characteristic speed $u$ does not have a definite sign. From
the theory of scalar conservation laws, it is seen that the boundary conditions  for $u, \rho$,  has to be prescribed in a weak sense \cite{b1}. Further, additional conditions  on the solutions need to be imposed
for uniqueness, these are called entropy conditions on the solutions; see Lax \cite{la1}, Hopf  \cite{h1}. In \cite{j3}
the solution of the inviscid system \eqref{1.1}, with Riemann type initial and boundary conditions  was constructed as the limit of solution of \eqref{1.3}, as $\epsilon$ goes to 0. This gives the entropy solution with a weak form of boundary conditions. In the analysis in \cite{j3},  apart from the formation of $\delta$ shocks, there is an additional  phenomenon of the formation of the boundary layer near the boundary at $x = 0$, that was observed, in the passage to the inviscid limit. In this paper we study the general initial and boundary conditions.

One basic question is the formulation of the solution for such problems. In this connection, the first such formulation was given in Lefloch \cite{le1} and a different formulation was given in Joseph \cite{j2} for the one dimensional $2 \times 2$ system. The approach of Lefloch was for a slightly general case and uses Lax formula, whereas that of Joseph was for one dimensional case of \eqref{1.1} for the Riemann problem and uses modified adhesion model.  

The works of  Danilov and  Shelkovich \cite{ds1, ds2} gave a new direction  by giving new definitions of weak formulation of the solution in the framework of weak asymptotic method. There are different notions of solutions introduced in \cite{j2, le1, ds1,ds2} and they contain interconnections and applications to different systems. First we  explain these different notions.

\textbf{Definition 1.1:} A family of smooth functions $\left( u^{\epsilon},\rho^{\epsilon} \right)_{\epsilon >0}$ is called a \textit{weak asymptotic solution} to the system
\eqref{1.1} with initial conditions $u(x,0)=u_{0}(x) \ , \ \rho(x,0)=\rho_{0}(x)$, and boundary conditions
 $u(0,t)=u_{B}(\cdot) \ , \ \rho(0,t)=\rho_{B}(t)$ provided: \\
i) As $\epsilon \rightarrow 0$ we have, for all  $\phi \in C_c^{\infty}(0, \infty),$ 
\begin{equation}
\begin{aligned}
\int_0^{\infty}{ ( u^\epsilon_t(x, t) + u^\epsilon(x, t)  u^\epsilon_x(x, t) ) \phi(x) dx} &=o(\epsilon),\\
\int_0^{\infty}{ ( \rho^\epsilon_t(x, t)  + {(\rho^\epsilon(x, t)  u^\epsilon(x, t) )}_x ) \phi(x)  \,\,\, dx} &= o(\epsilon),\\  
\end{aligned}
\label{1.6}
\end{equation}
The estimates \eqref{1.6} are required to hold uniformly on $[0, T]$ for each $T$.\\
ii) As $\epsilon \rightarrow 0,$ for all  $ \psi \in C_c^{\infty}(0, \infty), $ we have, 
\begin{equation}
\begin{aligned}
\int_0^{\infty}{ ( u^\epsilon(x,0)-u_0(x) ) \psi(x)  dx} &=o(\epsilon),\\
\int_0^{\infty} {( \rho^\epsilon(x,0)-\rho_0(x) ) \psi(x)  dx} &=o(\epsilon),\\
\int_0^{\infty} {( u^\epsilon(0,t)-u_B(t)) \psi(t)  dt} &=o(\epsilon),\\
\int_0^{\infty}{ ( \rho^\epsilon(0,t)-\rho_B(t)) \psi(t)  dt } &=o(\epsilon).
\end{aligned}
\label{1.7}
\end{equation}
We now define the notion of  a  generalized solution to the initial boundary value  problem  as follows.

\textbf{Definition 1.2:} A distribution $(u, \rho)$ is called a {\textit{generalized solution}} to the initial boundary value  problem in the domain $\Omega_T =\{(x,t) : x >0, t \in [0, T]\}$ if it is the limit in distribution of a weakly asymptotic solution $(\rho^\epsilon,u^\epsilon)_{\epsilon > 0}$ to the  initial boundary value problem, as $\epsilon$ goes to zero. 

Our aim is to use modified adhesion approximations to construct weak asymptotic solutions
for general $L^\infty$ initial data. This we do in section 2. 
In section 3, we study  interaction of boundary Riemann problem originating at $(0,0)$ and another solution with Riemann data originating at $(x_0,0), x>0.$

\section {Weak asymptotic solution of the initial boundary value problem}

In this section, we construct weak asymptotic solution to \eqref{1.1} with initial condition \eqref{1.4} and boundary condition \eqref{1.5}. We prove the following result.
\begin{theorem}
Let $(u_0,\rho_0)$ and $(u_B,\rho_B)$ are bounded measurable functions on $[0,\infty)$ with $\rho_0$ is integrable. Then there exists a weak
asymptotic solution $(u^\epsilon,\rho^\epsilon)$ to the equation \eqref{1.3} with initial condition (1.4) and boundary condition (1.5), in the sense of  definition (1.1).
\end{theorem}
 
The proof of this result is based on existence and smoothness of solution to the modified 
adhesion approximation \eqref{1.3} with initial and 
boundary conditions \eqref{1.4} and \eqref{1.5} for smooth data,
$(u_0(x)$, $\rho_0(x))$ and $(u_B(t),\rho_{B}(t))$. More precisely, first we  prove the following theorem.
\begin{theorem}
Assume $u_0$ and $\rho_0$ are functions of bounded variation, $u_B$ is a 
constant and 
$\rho_B(t)$ is a Lipschitz continuous function, then the 
solution 
to the problem \eqref{1.3} with initial and boundary conditions \eqref{1.4} and \eqref{1.5} is given by 
\begin{equation}
\begin{aligned}
u =-\epsilon \frac{p_x}{p}, \,\,\,\rho=(\frac{q}{p})_x
\end{aligned}
\label{2.1}
\end{equation}
where $p$ and $q$ are given by the solutions of initial boundary value problems for heat equation
\begin{equation}
\begin{aligned}
&p_t=\frac{\epsilon}{2} p_{xx}, x>0, t>0\\
&p(x,0) = e^{-\frac{U_0(x)}{\epsilon}},,\,\,x>0,\\
&\epsilon p_x(0,t)+u_B(t)p(0,t)=0, t>0
\end{aligned}
\label{2.2}
\end{equation}
and
\begin{equation}
\begin{aligned}
&q_t =\frac{\epsilon}{2}q_{xx},\\
&q(x,0)= R_0(x) e^{-\frac{U_0(x)}{\epsilon}},\\ 
&\epsilon q_x(0,t)+u_B(t)q(0,t)=\epsilon p(0,t)\rho_B(t),
\end{aligned}
\label{2.3}
\end{equation}
where $R_0(x)$ and $U_0(x)$ are given by
\begin{equation}
 U_0(x)=\int_0^x u_0(s)ds,\,\,\,R_0(x)=\int_0^x \rho_0(s)ds
\label{2.4}
\end{equation}
\end{theorem}

 { \bf Proof:} First we take arbitrary bounded measurable 
functions, $u_0(x),\rho_0(x), u_B(t)$ and $\rho_B(t)$. If $U(x,t), R(x, t)$ are 
solutions to 
\begin{equation}
\begin{aligned} 
&U_t + \frac{(U_x)^2}{2} = \frac{\epsilon}{2}U_{xx},\\
&R_t + u R_x = \frac{\epsilon}{2}R_{xx} 
\label{2.5}
\end{aligned} 
\end{equation}
in $\{(x,t) :x>0,t>0\}$ with initial and boundary conditions
\begin{equation}
\begin{aligned} 
&U(x,0) = U_0(x),\,\,R(x,0)=R_0(x)\\
&U_x(0,t) = u_B(t),\,\,R_x(0,t)=\rho_B(t)
\label{2.6}
\end{aligned} 
\end{equation} 
then 
\begin{equation} 
(u(x,t),\rho(x,t) = (U_{x}(x,t),R_x(x,t))
\label{2.7} 
\end{equation} 
is a solution of \eqref{1.3} with initial and boundary data \eqref{1.4} \eqref{1.5}, where $R_0$ and $U_0$ are
given by \eqref{2.4}.
We introduce  new unknown variables $(p, q)$. The 
unknown $p$ is defined by the usual Hopf-Cole transformation and $q$, 
by a modified version of it in the following way 
\begin{equation} 
p = e^{-\frac{U}{\epsilon}}, 
q = R e^{-\frac{U}{\epsilon}}.
\label{2.8}
\end{equation}
An easy calculation shows that
\begin{equation}
p_t - \frac{\epsilon}{2}p_{xx} = -\frac{1}{\epsilon}[U_t +
\frac{(U_x)^2}{2}
- \frac{\epsilon}{2}U_{xx}] e^{-\frac{w}{\epsilon}},
\label{2.9}
\end{equation}
and 
\begin{equation} 
q_t - \frac{\epsilon}{2}q_{xx} = [R_t + u R_{x} -\frac{\epsilon}{2}
R_{xx}] e^{-\frac{w}{\epsilon}} 
- \frac{1}{\epsilon}[U_t + \frac{(U_x)^2}{2} - \frac{\epsilon}{2}U_{xx}] 
U e^{-\frac{U}{\epsilon}}.
\label{2.10}
\end{equation}
It follows from \eqref{2.8}-\eqref{2.10} 
that $U$ and $R$ are 
solutions of \eqref{2.5}-\eqref{2.6} iff $p$ 
and $q$ are solutions of linear problems
\begin{equation}
\begin{aligned}
&p_t=\frac{\epsilon}{2} p_{xx},\\
&p(x,0) = e^{-\frac{U_0(x)}{\epsilon}},\\
&\epsilon p_x(0,t)+u_B(t)p(0,t)=0
\end{aligned}
\label{2.11}
\end{equation}
and
\begin{equation}
\begin{aligned}
&q_t =\frac{\epsilon}{2}q_{xx},\\
&q(x,0)= R_0(x) e^{-\frac{U_0(x)}{\epsilon}},\\ 
&\epsilon q_x(0,t)+u_B(t)q(0,t)=\epsilon p(0,t)\rho_B(t),
\end{aligned}
\label{2.12}
\end{equation}
where we have also used \eqref{2.6} and \eqref{2.7}.

By  existence uniqueness and regularity theory of linear parabolic equations \cite{Friedman} there exists unique classical solution
$p,q$ to these initial boundary value problems. Further $p$ and $q$ are $C^\infty((0,\infty)\times(0,\infty))$ and $p>0$ in $x\geq 0, t\geq 0$. Since  $(u,\rho)$ can be expressed in terms of $(p,q)$ by the transformation 
\eqref{2.7} and \eqref{2.8} we have a $C^\infty$ solution $(u,\rho)$.  
This completes the proof of the theorem. \\

 We use this theorem to give a proof of Theorem 2.1.

 { \bf Proof of Theorem 2.1:  }\\
First we regularise the initial and boundary data data using a cut off near $0$ and then regularize by convolution in the scale $\epsilon$ as follows. Let $\chi_{[2\epsilon,\infty)}$ be the characteristic function on $[ 2 \epsilon, \infty)$ and $\eta_\epsilon$ is the usual Friedrichs mollifier in one space dimension. Let
\begin{equation}
\begin{aligned}
&u^\epsilon_0(x)=(u_0 \chi_{[2\epsilon,\infty)}*\eta^\epsilon)(x),\,\,\, \rho^\epsilon_0(x)=(\rho_0 \chi_{[2\epsilon,\infty)}*\eta^\epsilon)(x),\\
&u^\epsilon_B(t)=(u_B \chi_{[2\epsilon,\infty)}*\eta^\epsilon)(t),\,\,\, \rho^\epsilon_B(t)=(u_B \chi_{[2\epsilon,\infty)}*\eta^\epsilon)(t)
\end{aligned}
\label{2.13}
\end{equation}
Now let $(u^\epsilon,\rho^\epsilon)$ be the solution to \eqref{1.3} with initial conditions
\[
u^\epsilon(x,0)=u^\epsilon_0(x),\,\,\,\rho^\epsilon(x,0)=\rho^\epsilon_0(x), x>0
\]
and boundary conditions 
\[
u^\epsilon(0,t)=u^\epsilon_B(t),\,\,\,\rho^\epsilon(0,t)=\rho^\epsilon_B(t), t>0
\]
constructed by the Theorem 2.2. Then clearly $(u^\epsilon,\rho^\epsilon)$  is a $C^\infty$ function. By maximum principle $u^\epsilon$ is  bounded by $c_1 = \max \{||u_0||_{L^\infty(0,\infty)}, ||u_B||_{L^\infty(0,\infty)}\}$. Further  $\rho^\epsilon$ is the space derivatve of a smooth function $R^\epsilon$ which satisfies
\[
\begin{aligned} 
&R_t + u R_x = \frac{\epsilon}{2}R_{xx} , \,\, x>0,t>0\}\\
&R(x,0)=R^\epsilon_0(x),x>0\\\
&R_x(0,t)=\rho^\epsilon_B(t),t>0.
\end{aligned}
\]
We show that $R^\epsilon$ is bounded independent of $\epsilon$, on bounded subsets of $[0,\infty) \times (0,\infty)$ using
standard comparison theorems, see \cite{Lieb}. To do this we apply the partial differential operator 
$L=\epsilon \partial_{xx} - u^\epsilon \partial_{x} -\partial_{t}$  on the function  $w(x,t) = c_2 (x+1) +c_3(t+1)$, with 
$c_2=||R_0||_{L^\infty(0,\infty)}+||\rho_B||_{L^\infty(0,\infty)}$ and $c_3>c_1c_2 $, we get $L w=-c_2 u^\epsilon -c_3 \leq 0 = LR^\epsilon$ and  $w_x(0,t) =c_2>\rho_B(t)=R_x^\epsilon( 0,t)$ and $w(x,0)=c_2(x+1) \geq R_0(x)=R^\epsilon(x,0)$ . By comparison theorem we get $R^\epsilon(x,t) \leq w(x,t)$. Similarly working with the comparison function $-w(x,t)$ we get 
$R^\epsilon(x,t) \geq -w(x,t)$. Combining the two we get 
\[
-  (c_2 (x+1) +c_3(t+1)) \leq R^\epsilon(x,t)\leq  c_2 (x+1) +c_3(t+1).
\]

 Now to show $(u^\epsilon,\rho^\epsilon)$ is a weak asymptotic solution of the equation (1.1), it is enough to note that
for $\phi$ any  $C^\infty$ function with compact support in $(0,\infty)$
\[
\begin{aligned}
\int_0^{\infty}(u^\epsilon_t(x, t) + u^\epsilon(x, t)  u^\epsilon_x(x, t)) \phi(x) dx& =\epsilon \int_0^\infty u^\epsilon_{xx}(x,t)\phi(x) dx\\
& =\epsilon \int_0^\infty u^\epsilon(x,t)\phi_{xx}(x) dx\\
\int_0^{\infty}( \rho^\epsilon_t(x, t)  + {(\rho^\epsilon(x, t)  u^\epsilon(x, t) )}_x  \phi(x))dx&=\epsilon \int_0^\infty \rho^\epsilon_{xx}(x,t)\phi(x) dx \\
&=-\epsilon \int_0^\infty R^\epsilon(x,t)\phi_{xxx}(x) dx.
\end{aligned}
\]
Since $u^\epsilon$ is bounded and $R^\epsilon$ is bounded on compact subsets of $[0,\infty) \times [0,\infty)$ the righthand sides  goes to zero uniformly in $t \in [0,T]$, as $\epsilon$ tends to $0$, for all  $C^\infty$ test functions $\phi$,  supported  in $(0,\infty)$.
The initial and boundary conditions satisfies in the sense of Definition 1.1, now easily follows from the regularization property by convolution.\\

{\bf Remark :} For special initial data we can say more. When the boundary data $u_B$ is constant, the solution $(u^\epsilon,\rho^\epsilon)$ given in theorem can be explicitly written down, see \cite{Andrei}:
 \begin{equation}
\begin{aligned}
&p^\epsilon(x,t)=\int_0^\infty e^{-\frac{U_0(y)}{\epsilon}} 
K(x,y,t,\epsilon) 
dy\\
&q^\epsilon(x,t)=\int_0^\infty R_0(y) 
e^{-\frac{U_0(y)}{\epsilon}}K(x,y,t,\epsilon) dy
-\frac{\epsilon}{2} \int_0^t K(x,0,t-\tau,\epsilon) 
p^\epsilon(0,\tau)\rho_B(\tau) d\tau
\end{aligned}
\label{2.14}
\end{equation}
with
\begin{equation}
R_0(x) =\int_0^x \rho_0(y) dy,\,\,\,U_0(x)=\int_0^x u_0(y) dy,
\label{2.15}
\end{equation}
and $K(x,y,t,\epsilon)$ given by ,

\begin{equation}
\begin{aligned}
K(x,y,t,\epsilon)&=\frac{1}{\sqrt(2 \pi t 
\epsilon)}(e^{-\frac{(x-y)^2}{2t\epsilon}} 
+e^{-\frac{(x+y)^2}{2t\epsilon}})+
 \frac{2u_B/\epsilon}{\sqrt(2 \pi t 
\epsilon)}\int_0^\infty 
e^{-\frac{1}{\epsilon}\{\frac{(x+y+z)^2}{2 t}-u_B z\}} dz.
\end{aligned}
\label{2.16}
\end{equation} 
The vanishing viscosity solution was proved  for the Riemann type initial data   in \cite{j3} and given by the following theorem. Depending upon the speed there are  three types of waves, classical shocks, rarefaction or delta waves. 
\begin{theorem}
Let $u^\epsilon$ and $\rho^\epsilon$ be the solution of 
\eqref{1.3} with initial condition,
\[
(u(x,0),\rho(x,0)) = (u_0,\rho_0), \,\,\, x>0,
\]
 and boundary condition
\[
(u(0,t),\rho(0,t) = (u_b,\rho_b), \,\,\, t>0,
\]
 with $u_0,\rho_0,u_b,\rho_b$ are
all constants then $u(x,t)= \lim_{\epsilon \rightarrow 0} u^\epsilon(x,t)$
exists pointwise a.e. and $\rho(x,t)= lim_{\epsilon \rightarrow 0} 
\rho^\epsilon(x,t)$, in the sense of distributions and $(u,\rho)$ 
have the following form:

{Case 1: (Refer Figure 1.)  $u_0=u_b>0$.} Here characteristics are all straight lines with speed $u_0$ and
\[
(u(x,t),\rho(x,t)) = \begin{cases} \displaystyle
      {(u_0,\rho_b) ,\,\,\,if \,\,\,x < u_0 t} 
\\\displaystyle
        {(u_0, \rho_0),\,\,\,if \,\,\, x > u_0 t}.
\end{cases}
\]

{Case 2: $u_0\leq 0,u_b \leq0$,}
\[
(u(x,t),\rho(x,t) = (u_0,\rho_0).
\]

{Case 3: (Refer Figure 2.) $0<u_b<u_0$,} characteristics are all straight lines with speed $u_0$ and
\[
(u(x,t),\rho(x,t)) = \begin{cases} \displaystyle
      {(u_b,\rho_b),\,\,\, if \,\,\,x<u_b t} 
\\\displaystyle
        { (x/t,0),\,\,\, if \,\,\,u_bt< x<u_0 t}
\\\displaystyle
        {(u_0,\rho_0),\,\,\, if \,\,\, x>u_0 t}.
\end{cases}
\]

\begin{figure}[!hbtp]
\includegraphics[width=12cm,height=5cm]{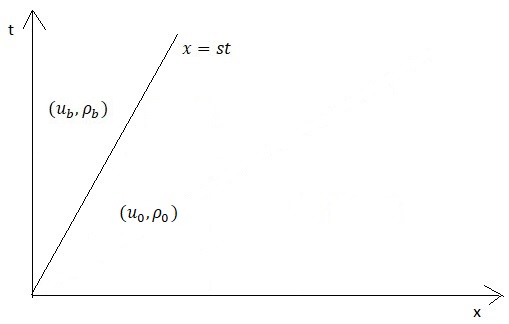}
\caption{ Shock wave}
\end{figure} 

{Case 4:  $u_b<0<u_0$,} 
\[
(u(x,t),\rho(x,t)) = \begin{cases} \displaystyle
        {(x/t, 0),\,\,\, if \,\,\,0< x<u_0 t}
\\\displaystyle
        {(u_0,\rho_0),\,\,\, if \,\,\,x>u_0 t}.
\end{cases}
\]

{Case 5: $u_0<u_b$ and $u_b+u_0>0$, }
\[
(u(x,t),\rho(x,t)) = \begin{cases} \displaystyle
      {(u_B,\rho_b),\,\,\, if \,\,\,x<st} 
\\\displaystyle
        {(1/2(u_b+u_0), (1/2)(u_b-u_0)(\rho_0+\rho_b) t \delta_{x=st})
                      ,\,\,\, if \,\,\,x=st}
\\\displaystyle
        { (u_0,\rho_0),\,\,\, if \,\,\,x>st}.
\end{cases}
\]

\begin{figure}[!hbtp]
\includegraphics[width=12cm,height=4cm]{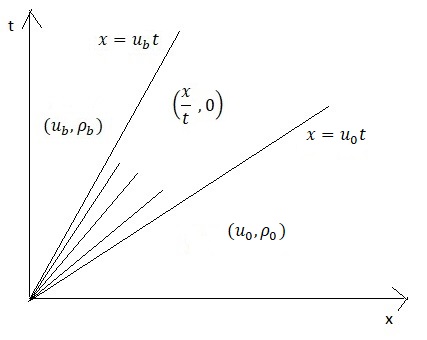}
\caption{ Rarefaction wave}
\end{figure}

where $s=\frac{u_0+u_b}{2}$. This limit satisfies the equation 
\eqref{1.1} and the initial condition \eqref{1.4} and a weak form of the boundary condition \eqref{1.5} namely,
\begin{equation}
u(0+,t) \in E(u_b),\rho(0+,t)) = \rho_b,\,\,\, if \,\,\, u(0+,t)>0
\label{e2.17}
\end{equation}
where for each $u_b \in R^1$, the admissible set $E(u_b)$ is defined by
\begin{equation}
E(u_b) =  \begin{cases} \displaystyle
      {(-\infty,0]) ,\,\,\,if \,\,\,u_b \leq 0}
\\\displaystyle
        {\{u_b\}\cup (-\infty, -u_b],\,\,\,if \,\,\, u_b>0}.
\end{cases}
\label{e2.18}
\end{equation}
\end{theorem}

\section{Interaction of Riemann and Riemann Boundary value problem}
In this section we would like to study interaction of boundary value problem in $x>0,t>0$ with a Riemann problem originating at a point $(x_0,0), x_0>0$. More pecisely, we construct  the solution to  \eqref{1.1} with initial condition
\begin{equation}
(u(x,0),\rho(x,0)) = \begin{cases} \displaystyle
      {(u_L,\rho_L) ,\,\,\,if \,\,\,0<x < x_0} 
\\\displaystyle
        {(u_R, \rho_R),\,\,\,if \,\,\, x > x_0 }
\end{cases}
\label{e3.1}
\end{equation}
and boundary condition
\begin{equation}
(u(0,t),\rho(0,t)) =  \displaystyle
      {(u_b,\rho_b) ,\,\,\,if \,\,\,t>0} 
\label{e3.2}
\end{equation}
 Here the solution is a result of interaction of waves originationg  at $(0,0)$  and at $(x_0,0)$, $x_0>0$.  We use   theorem 2.2, and  the 
following theorem which gives the solution to Riemann problem is proved in \cite{j1} with initial discontinuity at $x=x_0$, $x_0>0$ by shifting the origin to $x_0$ by a change of variable $x$ to $x-x_0$.
\begin{theorem}
Let $u^\epsilon$ and $\rho^\epsilon$ be the solution of
\eqref{1.3} with initial condition,
\[
(u(x,0),\rho(x,0)) = \begin{cases} \displaystyle
      {(u_L,\rho_L) ,\,\,\,if \,\,\,x < 0} 
\\\displaystyle
        {(u_R, \rho_R),\,\,\,if \,\,\, x > 0 }.
\end{cases}
\]
 with $x_0>0$ and $u_L,\rho_L,u_R,\rho_R$ are
all constants then $u(x,t)= \lim_{\epsilon \rightarrow 0} u^\epsilon(x,t)$
exists pointwise a.e. and $\rho(x,t)= lim_{\epsilon \rightarrow 0} 
\rho^\epsilon(x,t)$, in the sense of distributions and $(u,\rho)$ 
have the following form:

{Case 1: $u_L=u_R=u_0 >0$,} 
\[
(u(x,t),\rho(x,t)) = \begin{cases} \displaystyle
      {(u_0,\rho_L) ,\,\,\,if \,\,\,x < u_0 t} 
\\\displaystyle
        {(u_0, \rho_R),\,\,\,if \,\,\, x > u_0 t}.
\end{cases}
\]
{Case 2: $u_L<u_0$,}
\[
(u(x,t),\rho(x,t)) = \begin{cases} \displaystyle
      {(u_L,\rho_L),\,\,\, if \,\,\,x<u_L t} 
\\\displaystyle
        { (x/t,0),\,\,\, if \,\,\,u_L t< x<u_R t}
\\\displaystyle
        {(u_R,\rho_R),\,\,\, if \,\,\, x>u_R t}.
\end{cases}
\]

{Case 3: $u_R<u_L$}
\[
(u(x,t),\rho(x,t)) = \begin{cases} \displaystyle
      {(u_L,\rho_L),\,\,\, if \,\,\,x<st} 
\\\displaystyle
        {(1/2(u_L+u_R), (1/2)(u_L-u_R)(\rho_L+\rho_R) t \delta_{x=st})
                      ,\,\,\, if \,\,\,x=st}
\\\displaystyle
        { (u_R,\rho_R),\,\,\, if \,\,\,x>st},
\end{cases}
\]
where $s=\frac{u_L+u_R}{2}$. This limit satisfies the equation 
\eqref{1.1} and the initial condition \eqref{1.4}.
\end{theorem}

We study interaction by piecing together  these two solutions in the previous section. Here we use the following  result,  from
\cite{ds1}.
Let $(u,\rho)$ be a distribution which is smooth except along a curve $x=s(t)$ in some neighborhood of a point $(s(t_0),t_0)$ and having the form 
\begin{equation*}
\begin{aligned}
u (x,t) &= u_{l} (x,t) +H(x-s(t)) ( u_{r} (x,t))- u_{l} (x,t) ), \\
\rho (x,t) &= \rho_{l}(x,t) +H(x-s(t)) (  \rho_{r}(x,t) - \rho_{l}(x,t)) +e(t) \delta_{x=s(t)}.
\end{aligned}
\end{equation*}
where $H$ denotes the Heaviside function. Then $(u,\rho)$ is a solution to \eqref{1.1} in that neighborhood   if and only if $(u,\rho)$ is a classical solution in the neighborhood of $(s(t_0),t_0)$ except along $x=s(t)$ and along $x=s(t)$, $(u,\rho)$ satisfies the Rankine-Hugoniot conditions 

\begin{equation}
\begin{aligned}
s'(t)&=\frac{u(s(t)+,t) +u(s(t)-,t)}{2},\\
e'(t)&= -s'(t)(\rho(s(t)+,t)-\rho(s(t)-,t))  +(\rho u)(s(t)+,t)-(\rho u)(s(t)-,t).
\end{aligned}
\label{e3.3}
\end{equation}

\begin{theorem}
Solution to the system \eqref{1.1} with initial condition \eqref{e3.1} and boundary conditions \eqref{e3.2} exists and is given by \eqref{e3.4} - \eqref{e3.31}. These solutions has the following properties. 
\begin{itemize}
\item When the speed $u$ on the left is more than the speed on the right , the particles collide, shocks are formed and there is concentration of mass on the shock curve due to which  density $\rho$ becomes a measure concentrated on this shock  curve, called delta waves. The speed of the shock $s'(t)$ and the strength of the delta wave $e(t)$ satisfy \eqref{e3.3}.
\item In the rarefaction region the density $\rho$ is always zero.
 \item Three types of waves are possible originating at $(0,0)$ and $(x_0,0)$ namely, shock wave, rarefaction wave and delta wave, depending on initial and boundary values of the speed $u$. 
\item  The soutions satisfy boundary conditions in the weak sense \eqref{e2.17}-\eqref{e2.18}. When $u_b\leq0$, there is no effect of boundary values $(u_b,\rho_b)$  and the solution  is same as the case  $u_b=0$.
\item The interaction of  these waves happens when the speed of the wave on the left is more than the speed on the right. 
\item When a delta wave interacts with a rarefaction or another delta wave the path and strength of the delta wave changes but it can never disappear in the region $x>0,t>0$, if initial and boundary densites are positive.

\end{itemize}
\end{theorem}
{ \bf Proof:} 
Case 1: $u_L = u_b >0:$ \\
Subcase 1: $u_L = u_b= u_R=  u_0>0.$ (Refer Figure 3.) In this case the solution behaves as linear equation and \\
the solution in this case is, 
\begin{equation}
(u(x,t),\rho(x,t)) = \begin{cases} \displaystyle
      {(u_0,\rho_b),\,\,\, if \,\,\,x<u_0 t, \,\,\, \forall t, } 
\\\displaystyle
        { (u_0,\rho_L),\,\,\, if \,\,\,u_0 t< x<u_0 t +x_0,  \,\,\, \forall t,}
\\\displaystyle
        {(u_0,\rho_R),\,\,\, if \,\,\, x>u_0 t +x_0,  \,\,\, \forall t.}
\end{cases}
\label{e3.4}
\end{equation}

\begin{figure}[!hbtp]
\includegraphics[width=13cm,height=5cm]{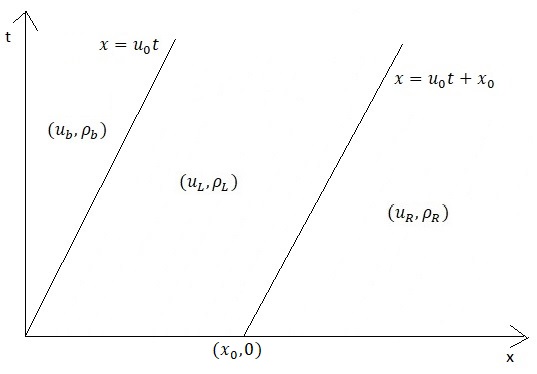}
\caption{ Constant speed $u_0$, no interaction}
\end{figure}

Subcase 2: $u_L = u_b=u_0>0, \,\,\,  0< u_L < u_R.$ In this case the speed of characteristics starting on the left of $x_0$ 
is less than those on the right and hence
$(u_L, \rho_L)$ is connected to $(u_R, \rho_R)$ by a rarefaction centered at $(x_0,0)$. The solution  is, 
\begin{equation}
\begin{aligned}
( u (x, t), \rho(x, t))= \begin{cases}
(u_0, \rho_b), &\text{ if } 0 < x < u_L t, \,\,\,  \forall t, \\
(u_0, \rho_L), &\text{ if }  u_L t < x < u_L t + x_0, \,\,\,  \forall t, \\
( \frac{x-x_0}{t}, 0), &\text{ if }   u_L t  + x_0< x < u_R t + x_0,\,\,\,   \forall t, \\
(u_R, \rho_R),  &\text{ if } x> u_R t + x_0,\,\,\,   \forall t.
\end{cases}
\label{e3.5}
\end{aligned}
\end{equation}

Subcase 3:  $u_L = u_b=u_0>0, \,\,\,  u_L > u_R.$ Here characteristic starting on the left of $(x_0,0)$ meet those starting  on the right and a shock is formed starting at $(x_0,0)$ namely  $x =  (\frac{ u_L + u_R }{2}) t + x_0$. The solution contains a delta measure concentrated on the the shock curve.  Let $(x_1,t_1)$ be the meeting point of $x+u_b t$ and $x-\frac{u_L+u_R}{2} t+x_0$, then $(x_1,t_1)=(\frac{2u_L x_0}{u_L-u_R},(\frac{2 x_0}{u_L-u_R})$. The solution is
\begin{equation*}
( u (x, t), \rho(x, t)) 
\end{equation*}
\begin{equation}
\begin{aligned}
 = \begin{cases}
(u_0, \rho_b), &\text{ if } x< u_0 t,\\
(u_0, \rho_L), &\text{ if }  u_0 t < x < (\frac{ u_0 + u_R }{2}) t + x_0, \\
( \frac{1}{2} (u_L+u_R), \frac{1}{2} (u_L- u_R)(\rho_L+\rho_R) t \delta_{x= (\frac{ u_L + u_R }{2}) t}),   &\text{ if }  x = (\frac{ u_L + u_R }{2}) t + x_0, x>0,\\
(u_R, \rho_R), &\text{ if }   x > (\frac{ u_L + u_R }{2}) t + x_0,\\
(u_0, \rho_b), &\text{ if } x<  (\frac{ u_0 + u_R }{2}) t + x_0, \,\,\, t > t_1,\\
( \frac{1}{2} (u_0+u_R), \frac{1}{2} (u_0 - u_R)(\rho_b+\rho_R) t \delta_{x= (\frac{ u_L + u_R }{2}) t}),   &\text{ if }  x = (\frac{ u_L + u_R }{2}) t + x_0,\,\,\, t > t_1,\\
\end{cases}
\label{e3.6}
\end{aligned}
\end{equation}

Case 2: If $u_L = u_b <0:$ \\
Subcase 1: $u_L = u_b= u_R=u_0 <0.$ (Refer Figure 4.) This is as linear case and there is no interaction.
The solution in this case is, 
\begin{equation}
\begin{aligned}
( u (x, t), \rho(x, t))= \begin{cases}
(u_0, \rho_L), &\text{ if } 0<x< u_0 t +x_0, \\
(u_0, \rho_R),  &\text{ if } x> u_0 t + x_0.
\end{cases}
\label{e3.7}
\end{aligned}
\end{equation}
.

\begin{figure}[!hbtp]
\includegraphics[width=10cm,height=5cm]{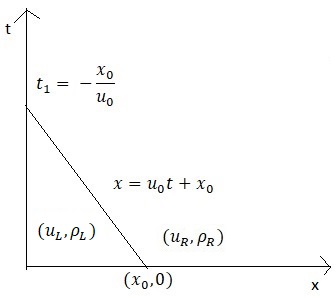}
\caption{Constant speed $u_0$, no interaction }
\end{figure}

Subcase 2: $u_L = u_b=u_0<0, \,\,\,   u_L < u_R.$ ( Refer Figure 5.) Here
$(u_L, \rho_L)$ is connected to $(u_R, \rho_R)$ by a rarefaction and the solution in this case is, 
\begin{equation}
\begin{aligned}
( u (x, t), \rho(x, t))= \begin{cases}
(u_0, \rho_L), &\text{ if } 0<x< u_L t, \\
( \frac{x-x_0}{t}, 0), &\text{ if }  u_L t +x_0 < x < u_R t + x_0, \\
(u_R, \rho_R),  &\text{ if } x> u_R t + x_0.
\end{cases}
\label{e3.8}
\end{aligned}
\end{equation}

\begin{figure}[!hbtp]
\includegraphics[width=12cm,height=7cm]{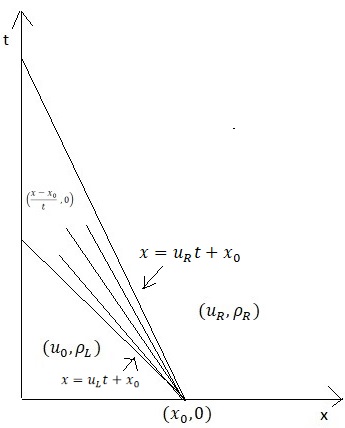}
\caption{ Rarefaction centered at $(x_0,0)$}
\end{figure}

Subcase 3: $u_L = u_b=u_0<0, \,\,\,  u_L > u_R.$ In this case 
$(u_L, \rho_L)$ is connected to $(u_R, \rho_R)$ by a shock with delta measure concentrated on this shock curve. The solution is, 
\begin{equation*}
( u (x, t), \rho(x, t)) 
\end{equation*}
\begin{equation}
\begin{aligned}
= \begin{cases}
(u_0, \rho_L), &\text{ if }  0 < x < (\frac{ u_L + u_R }{2}) t + x_0,  \\
( \frac{1}{2} (u_L+u_R), \frac{1}{2} (u_L- u_R)(\rho_L+\rho_R) t \delta_{x= (\frac{ u_L + u_R }{2}) t}),   &\text{ if }  x = (\frac{ u_L + u_R }{2}) t + x_0, \\
(u_R, \rho_R), &\text{ if }   x > (\frac{ u_L + u_R }{2}) t + x_0.
\end{cases}
\end{aligned}
\label{e3.9}
\end{equation}

Case 3: If $u_L = u_R >0:$ \\
Subcase 1: If $u_L = u_R=u_0, \,\,\, u_L < u_b, u_L+u_b>0$ (Refer Figure 6.) Here a shock is formed at $(0,0)$ with shock curve
$x= \frac{ u_L + u_b}{2} t$ and this curve intersect with $x= u_L t + x_0$  at say $(x_1, t_1).$ For $t<t_1,$ $(u_b, \rho_b)$ is connected to $(u_L, \rho_L)$ by a shock. For $t>t_1, \,\,\, (u_b, \rho_b)$ is connected to $(u_R, \rho_R)$ by a shock and the solution in this case is, 
\begin{equation*}
( u (x, t), \rho(x, t)) 
\end{equation*}
\begin{equation}
\begin{aligned}
= \begin{cases}
(u_b, \rho_b), &\text{ if } x< (\frac{ u_L + u_b}{2}) t, \,\,\, 0< t< t_1, \\
( \frac{1}{2} (u_L+u_b), \frac{1}{2} (u_b- u_L)(\rho_b+\rho_L) t \delta_{x= (\frac{ u_L + u_b }{2}) t}),   &\text{ if }  x = (\frac{ u_L + u_b }{2}) t,  \,\,\, t< t_1, \\
(u_0, \rho_L), &\text{ if } (\frac{ u_L + u_b}{2}) t < x < u_0 t + x_0,    t< t_1, \\
(u_0, \rho_R),  &\text{ if } x> u_0 t + x_0, \,\,\,  t< t_1, \\
(u_b, \rho_b), &\text{ if } 0< x<( \frac{ u_R + u_b}{2}) t+ x_1,    t> t_1, \\
( \frac{1}{2} (u_R+u_b), \frac{1}{2} (u_b- u_R)(\rho_b+\rho_R) t \delta_{x= (\frac{ u_R + u_b }{2}) t}),   &\text{ if }  x = (\frac{ u_R + u_b }{2}) t,  \,\,\, t> t_1, \\
(u_R, \rho_R),  &\text{ if } x>  (\frac{ u_R + u_b}{2}) t+ x_1, \,\,\,   t> t_1.
\end{cases}
\end{aligned}
\label{e3.10}
\end{equation}

\newpage

\begin{figure}[!hbtp]
\includegraphics[width=14cm,height=5cm]{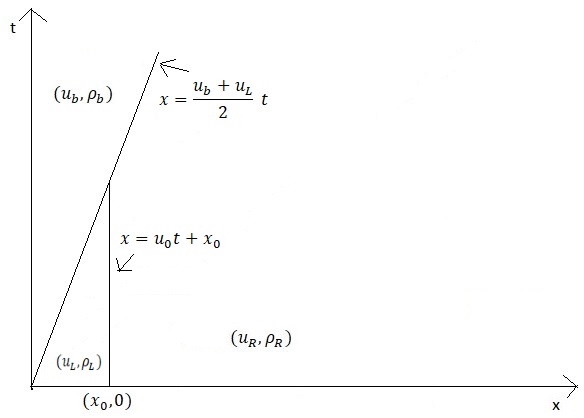}
\caption{ Shock wave, no interaction of waves}
\end{figure} 

Subcase 2: (Refer Figure 7.) If $u_L = u_R=u_0 >0, \,\,\, u_L > u_b, \,\,\, u_b>0.$ Here
$(u_b, \rho_b)$ is connected to $(u_L, \rho_L)$ by a rarefaction and the solution in this case is, 
\begin{equation}
\begin{aligned}
( u (x, t), \rho(x, t))= \begin{cases}
(u_b, \rho_b), &\text{ if } x< u_b t,  \,\,\,  \forall t, \\
( \frac{x}{t}, 0), &\text{ if }  u_b t < x < u_L t, \,\,\,   \forall t, \\
 (u_0, \rho_L), &\text{ if }    u_0 t< x<  u_0 t + x_0, \,\,\,   \forall t, \\
(u_0, \rho_R),  &\text{ if }   x> u_0 t + x_0,\,\,\,   \forall t. \\
\end{cases}
\end{aligned}
\label{e3.11}
\end{equation}

\begin{figure}[!hbtp]
\includegraphics[width=14cm,height=5cm]{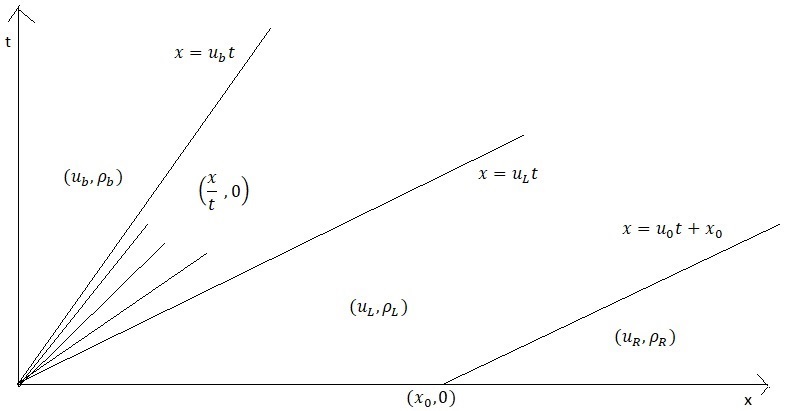}
\caption{ Rarefaction centered at $(0,0)$, no interaction}
\end{figure}

Subcase 3:  If $u_L = u_R=u_0 >0, \,\,\, u_L > u_b, \,\,\, u_b<0.$ (Refer Figure 7, with $u_b=0$).
The solution in this case is, 
\begin{equation}
\begin{aligned}
( u (x, t), \rho(x, t))= \begin{cases}
( \frac{x}{t}, 0), &\text{ if }  0 < x < u_0 t, \,\,\, \forall t, \\
 (u_0, \rho_L), &\text{ if }    u_0 t< x<  u_0 t + x_0 ,\,\,\,   \forall t, \\
(u_0, \rho_R),  &\text{ if }   x> u_0 t + x_0,\,\,\,   \forall t.
\end{cases}
\end{aligned}
\label{e3.12}
\end{equation}

Subcase 4:  If $u_L = u_R=u_0 <0, \,\,\, u_L < u_b, \,\,\, u_L+ u_b<0$ 
or if $u_L = u_R=u_0 <0, \,\,\, u_L > u_b, \,\,\, u_b < 0< u_L.$ 
The solution in this case is, 
\begin{equation}
\begin{aligned}
( u (x, t), \rho(x, t))= \begin{cases}
 (u_0, \rho_L), &\text{ if }    0< x<  u_0 t + x_0, \\
(u_0, \rho_R),  &\text{ if }   x> u_0 t + x_0  .
\end{cases}
\end{aligned}
\label{e3.13}
\end{equation}

Case 4: If $u_R= u_b >0:$ \\
Subcase 1:  If $u_R= u_b=u_0 >0, \,\,\, u_R < u_L.$ (Refer Figure 8.) In this case a rarefaction originating at $(0,0)$ and a shock wave origination at $(x_0,0)$ interact. This interaction starts at the point $(x_1, t_1)$ the point where 
the characeristic curve $x =u_L t$ and shock curve $x =  \frac{ u_R + u_L}{2} t +x_0$ intersect which can be explicitly calculated. The shock curve for $t > t_1$ is given by 
\begin{equation*}
\frac{dx}{dt}= \frac{1}{2} ( \frac{x}{t} + u_R), \,\,\, x(t_1)= x_1,
\end{equation*}
that gives the curve $x=\beta_1(t)$ i.e.
$x= \frac{u_R}{2} t + c {t}^{ \frac{1}{2}}, \,\,\, c =  x_1 {t_1}^{-\frac{1}{2}} -  \frac{u_R}{2} {t_1}^{ \frac{1}{2}}$. 
In $t<t_1, \,\,\, (u_b, \rho_b)$ is connected to $(u_L, \rho_L)$ by a rarefaction and $(u_L, \rho_L)$ is connected to $(u_R, \rho_R)$ by a shock with a delta wave. In $t>t_1$ the delta wave continues along the shock  curves  $x=\beta_1(t)$ with strength $e(t)=\frac{u_R \rho_R}{2}(t-t_1)+c \rho_R (t^{1/2}-t_1^{1/2}+\frac{(u_b-u_L)}{2}(\rho_L+\rho_b)t_1$ which is obtained from solving the equation for $e(t)$ in \eqref{e3.3} for $t>t_1$ with initial condition $e(t_1)=\frac{(u_b-u_L)}{2}(\rho_L+\rho_b)t_1$. The solution  is, 
\begin{equation*}
( u (x, t), \rho(x, t))
\end{equation*}
\begin{equation}
\begin{aligned}
= \begin{cases}
(u_b, \rho_b), &\text{ if }0< x<  u_b t, \forall t,  \\
( \frac{x}{t}, 0),  &\text{ if } u_b t< x< u_L t, \,\,\,  t<t_1 \\
(u_L, \rho_L),  &\text{ if }  u_L t< x< (\frac{ u_R + u_L}{2}) t +x_0,   t<t_1 \\
( \frac{1}{2} (u_L+u_R), \frac{1}{2} (u_L- u_R)(\rho_L+\rho_R) t \delta_{x= (\frac{ u_L + u_R }{2}) t}),   &\text{ if }   x = (\frac{ u_R + u_L}{2}) t +x_0,  \,\,\, t<t_1 \\
(u_R, \rho_R), &\text{ if } x> (\frac{ u_R + u_L}{2}) t + x_0, \,\,\, t<t_1\\
( \frac{x}{t}, 0),  &\text{ if } u_b t< x< \beta_1(t) , \,\,\,  t>t_1 \\
e(t) \delta_{x=\beta_1(t)},  &\text{ if } t>t_1,\\
(u_R,\rho_R), &\text{ if } x>\beta_1(t),\,\,\,t>t_1.
\end{cases}
\end{aligned}
\label{e3.14}
\end{equation}

\begin{figure}[!hbtp]
\includegraphics[width=12cm,height=7cm]{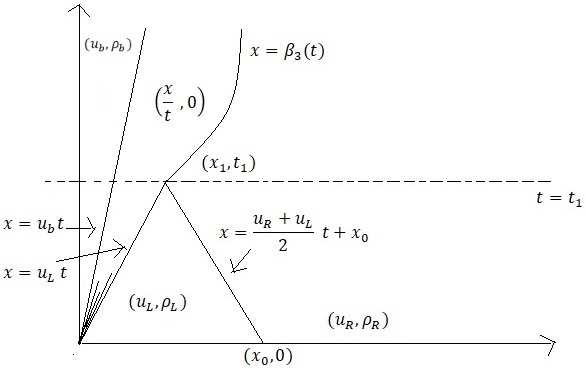}
\caption{ Interaction of rarefaction centered at $(0,0)$ and shock centered at $(x_0,0)$}
\end{figure}

Subcase 2: If $u_R= u_b=u_0 >0, \,\,\, u_R > u_L, \,\,\, u_b + u_L > 0.$ (Refer Figure 9.) Here the interaction is between shock originating at $(0,0)$  and rarefaction originating at $(x_0,0)$
The curves $x =u_L t$ and shock curve $x =  \frac{ u_b + u_L}{2} t $ intersect at say $(x_1, t_1)$. The curve for $t > t_1$ is given by 
\begin{equation*}
\frac{dx}{dt}= \frac{1}{2} ( \frac{x- x_0}{t} + u_b), \,\,\, x(t_1)= x_1,
\end{equation*}
that gives the curve $x=\beta_2(t)$ i.e.
$x= \frac{u_b}{2} t + c {t}^{ \frac{1}{2}}, \,\,\, c =  x_1 {t_1}^{-\frac{1}{2}} -  \frac{u_b}{2} {t_1}^{ \frac{1}{2}}$. In $t<t_1, \,\,\, (u_b, \rho_b)$ is connected to $(u_L, \rho_L)$ by a shock and $(u_L, \rho_L)$ is connected to $(u_R, \rho_R)$ by a rarefaction. In $t>t_1$ the shock curve continues as $x=\beta_2$ with shock strength  $e(t)=a t +b t^{1/2}+c$ where $a,b, c$ are constsnts depending only on initial data by solving \eqref{e3.3} as in previous case. The solution in this case is,  
\begin{equation*}
( u (x, t), \rho(x, t))
\end{equation*}
\begin{equation}
\begin{aligned}
= \begin{cases}
(u_b, \rho_b), &\text{ if }0< x<(  \frac{ u_b + u_L}{2}) t,   \,\,\, t<t_1 \\
( \frac{1}{2} (u_L+u_b), \frac{1}{2} (u_b- u_L)(\rho_L+\rho_b) t \delta_{x= (\frac{ u_L + u_b }{2}) t}),   &\text{ if }   x = (\frac{ u_b + u_L}{2}) t,   \,\,\, t<t_1 \\
(u_L, \rho_L),  &\text{ if } (\frac{ u_b + u_L}{2}) t< x< u_L t+ x_0,     t<t_1 \\
( \frac{x - x_0}{t}, 0),  &\text{ if }   u_L t+ x_0 < x<  u_R t+ x_0,   t<t_1 \\
(u_R, \rho_R), &\text{ if } x> u_R t+ x_0, \,\,\, t<t_1 \\
(u_b, \rho_b),  &\text{ if } 0< x< \beta_2(t), \,\,\,  t>t_1 \\
e(t) {\delta}_{x=\beta_2(t)},  &\text{ if } t > t_1, \,\,\, x= \beta_2(t), \\
( \frac{x}{t}, 0),  &\text{ if }   \beta_1(t)< x<   u_R t+ x_0,    t>t_1.
\end{cases}
\end{aligned}
\label{e3.15}
\end{equation}

\begin{figure}[!hbtp]
\includegraphics[width=15cm,height=9cm]{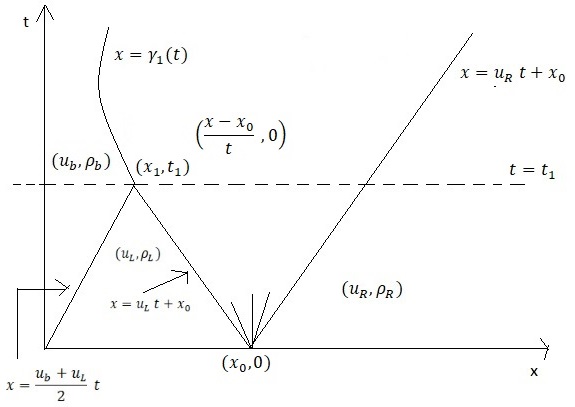}
\caption{ Interaction of shock centered at $(0,0)$ and rarefaction centered at $(x_0,0)$}
\end{figure}

Subcase 3:  If $u_R = u_b >0, \,\,\, u_R > u_L, \,\,\, u_L+ u_b<0.$ Here
$(u_L, \rho_L)$ is connected to $(u_R, \rho_R)$ by a rarefaction and  the solution in this case is, 
\begin{equation}
\begin{aligned}
( u (x, t), \rho(x, t))= \begin{cases}
 (u_L, \rho_L), &\text{ if }    0< x<  u_L t + x_0, \,\,\,  \forall t, \\
( \frac{x-x_0}{t}, 0),  &\text{ if }   u_L t+ x_0 < x<  u_R t+ x_0,   \,\,\, \forall t, \\
(u_R, \rho_R),  &\text{ if }   x> u_R t + x_0,\,\,\,  \forall t.
\end{cases}
\end{aligned}
\label{e3.16}
\end{equation}

Case 5: If $u_R= u_b <0:$ \\
Subcase 1: If $u_R= u_b <0, \,\,\, u_R < u_L, \,\,\, u_b < u_L.$ (Refer Figure 8 with $u_b=0$).Here rarefaction at $(0,0)$ interact with a shock at $(x_0,0)$. The starting point of interaction $(x_1, t_1)$ is the meeting point of
the curve $x =u_L t$ and shock curve $x =  \frac{ u_R + u_L}{2} t +x_0$ . The shock curve for $t > t_1$ is given by 
\begin{equation*}
\frac{dx}{dt}= \frac{1}{2} ( \frac{x}{t} + u_R), \,\,\, x(t_1)= x_1,
\end{equation*}
that gives the curve $x=\beta_1(t)$ i.e.
$x= \frac{u_R}{2} t + c {t}^{ \frac{1}{2}}, \,\,\, c =  x_1 {t_1}^{-\frac{1}{2}} -  \frac{u_R}{2} {t_1}^{ \frac{1}{2}}$.
As before   the solution in this case is, 
\begin{equation*}
( u (x, t), \rho(x, t))
\end{equation*} 
\begin{equation}
\begin{aligned}
= \begin{cases}
( \frac{x}{t}, 0),  &\text{ if } 0< x< u_L t, \,\,\,  t<t_1 \\
(u_L, \rho_L),  &\text{ if }  u_L t< x< (\frac{ u_R + u_L}{2}) t +x_0,  \\
( \frac{1}{2} (u_L+u_R), \frac{1}{2} (u_L- u_R)(\rho_L+\rho_R) t \delta_{x= (\frac{ u_L + u_R }{2}) t}),   &\text{ if }   x = (\frac{ u_R + u_L}{2}) t +x_0,  \\
(u_R, \rho_R), &\text{ if } x> (\frac{ u_R + u_L}{2}) t + x_0, \,\,\, t<t_1 \\
( \frac{x}{t}, 0),  &\text{ if } 0< x< \beta_3(t), \,\,\,  t>t_1 \\
e(t) {\delta}_{x=\beta_3(t)},  &\text{ if } t > t_1, \,\,\, x= \beta_3(t), \\
(u_R, \rho_R),  &\text{ if }  x>  \beta_3(t),  \,\,\,  t>t_1.
\end{cases}
\end{aligned}
\label{e3.17}
\end{equation}
where shock strength has the form $e(t)=a t +b t^{1/2}+c$ where $a,b, c$ are constants and obtained by solving \eqref{e3.3} with initial condition $e(t_1)= \frac{1}{2} (u_L- u_R)(\rho_L+\rho_R) t_1$.

Subcase 2:  If $u_R= u_b <0, \,\,\, u_R > u_L, \,\,\, u_b + u_L<0.$ Here
$(u_L, \rho_L)$ is connected to $(u_R, \rho_R)$ by a rarefaction and the solution in this case is, 
\begin{equation}
\begin{aligned}
( u (x, t), \rho(x, t))= \begin{cases}
 (u_L, \rho_L), &\text{ if }    0< x<  u_L t + x_0, \,\,\,  \forall t, \\
( \frac{x - x_0}{t}, 0),  &\text{ if }   u_L t+ x_0 < x<  u_R t+ x_0,   \,\,\, \forall t, \\
(u_R, \rho_R),  &\text{ if }   x> u_R t + x_0,\,\,\,  \forall t.
\end{cases}
\end{aligned}
\label{e3.18}
\end{equation}

Case 6: If $u_L \neq u_R \neq u_b:$ \\
Subcase 1: If $0< u_L < u_R < u_b$ 
or $ u_L <0<u_R< u_b, \,\,\, u_b +u_L >0$ 
or  $ u_L<u_R< 0 < u_b, \,\,\, u_b +u_L >0$. (Refer Figure 9.) In this case shock wave starting at $(0,0)$ interact with rarefaction wave at $(x_0,0)$.
The curves $x =u_L t + x_0$ and $x =  \frac{ u_b + u_L}{2} t +x_0$ intersect at say $(x_1, t_1)$. The curve for $t > t_1$ is given by 
\begin{equation*}
\frac{dx}{dt}= \frac{1}{2} ( \frac{x - x_0}{t} + u_b), \,\,\, x(t_1)= x_1,
\end{equation*}
that gives the curve $x=\gamma_1(t)$ i.e.
$x= \frac{u_b}{2} t + c {t}^{ \frac{1}{2}}, \,\,\, c =  x_1 {t_1}^{-\frac{1}{2}} -  \frac{u_b}{2} {t_1}^{ \frac{1}{2}}$.
In $t<t_1, \,\,\, (u_b, \rho_b)$ is connected to $(u_L, \rho_L)$ by a shock and $(u_L, \rho_L)$ is connected to $(u_R, \rho_R)$ by a rarefaction.   The solution in this case is, 
\begin{equation*}
( u (x, t), \rho(x, t))
\end{equation*}
\begin{equation}
\begin{aligned}
= \begin{cases}
(u_b, \rho_b), &\text{ if }0< x< ( \frac{ u_b + u_L}{2}) t,   \,\,\, t<t_1 \\
( \frac{1}{2} (u_L+u_b), \frac{1}{2} (u_b- u_L)(\rho_L+\rho_b) t \delta_{x= (\frac{ u_L + u_b }{2}) t}),   &\text{ if }   x = (\frac{ u_b + u_L}{2}) t,   \,\,\, t<t_1 \\
(u_L, \rho_L),  &\text{ if } (\frac{ u_b + u_L}{2}) t< x< u_L t+ x_0,   t<t_1 \\
( \frac{x-x_0}{t}, 0),  &\text{ if }   u_L t+ x_0 < x<  u_R t+ x_0,   t<t_1 \\ 
(u_R, \rho_R), &\text{ if } x> u_R t+ x_0, \,\,\,  \forall t \\
(u_b, \rho_b),  &\text{ if } 0< x< \gamma_1(t), \,\,\,  t>t_1 \\
e(t) {\delta}_{x=\gamma_1(t)},  &\text{ if } t > t_1, \,\,\, x= \gamma_1(t), \\
 ( \frac{x-x_0}{t}, 0),  &\text{ if }   \gamma_1(t)< x<   u_R t+ x_0,    t>t_1. 
\end{cases}
\end{aligned}
\label{e3.19}
\end{equation}
where shock strength has the form $e(t)=a t +b t^{1/2}+c$ where $a,b, c$ are constants and obtained by solving \eqref{e3.3} with initial condition $e(t_1)= \frac{1}{2} (u_b- u_L)(\rho_L+\rho_b) t_1$.

Subcase 2:  If $u_L< u_R< u_b <0$ 
or $ u_L<u_R<0< u_b, \,\,\, u_b +u_L <0$ 
or  $ u_L< 0<u_R< u_b, \,\,\, u_b +u_L <0.$ Here
$(u_L, \rho_L)$ is connected to $(u_R, \rho_R)$ by a rarefaction and the solution  is, 
\begin{equation}
\begin{aligned}
( u (x, t), \rho(x, t))= \begin{cases}
 (u_L, \rho_L), &\text{ if }    0< x<  u_L t + x_0, \\
( \frac{x - x_0}{t}, 0),  &\text{ if }   u_L t+ x_0 < x<  u_R t+ x_0,   \,\,\,  x>0, \\
(u_R, \rho_R),  &\text{ if }   x> u_R t + x_0,\,\,\,   \forall t.
\end{cases}
\end{aligned}
\label{e3.20}
\end{equation}

Subcase 3:  If $0< u_b< u_L< u_R.$ In this case,
$(u_b, \rho_b)$ is connected to $(u_L, \rho_L)$ by a rarefaction, $(u_L, \rho_L)$ is connected to $(u_R, \rho_R)$ by a rarefaction and the solution  is, 
\begin{equation}
\begin{aligned}
( u (x, t), \rho(x, t))= \begin{cases}
 (u_b, \rho_b), &\text{ if }    0< x<  u_b t, \,\,\, \forall t, \\
( \frac{x}{t}, 0),  &\text{ if }   u_b t< x<  u_L t,   \,\,\,   \forall t, \\
(u_L, \rho_L),  &\text{ if }   u_L t< x< u_L t + x_0 ,\,\,\,    \forall t, \\
( \frac{x- x_0}{t}, 0),  &\text{ if }   u_L t+ x_0 < x<  u_R t+ x_0, \,\,\, \forall t, \\
(u_R, \rho_R),  &\text{ if }   x> u_R t + x_0,\,\,\,   \forall t.
\end{cases}
\end{aligned}
\label{e3.21}
\end{equation}

Subcase 4:  If $ u_b<0< u_L< u_R.$ Here
 $(u_L, \rho_L)$ is connected to $(u_R, \rho_R)$ by a rarefaction and the solution is, 
\begin{equation}
\begin{aligned}
( u (x, t), \rho(x, t))= \begin{cases}
( \frac{x}{t}, 0),  &\text{ if }   0< x<  u_L t,   \,\,\,  \forall t, \\
(u_L, \rho_L),  &\text{ if }   u_L t< x< u_L t + x_0 ,\,\,\,   \forall t, \\
( \frac{x - x_0}{t}, 0),  &\text{ if }   u_L t+ x_0 < x<  u_R t+ x_0, \,\,\, \forall t, \\
(u_R, \rho_R),  &\text{ if }   x> u_R t + x_0,\,\,\,  \forall t.
\end{cases}
\end{aligned}
\label{e3.22}
\end{equation} 

Subcase 5:  If $ u_b< u_L<0< u_R$ 
or  $ u_b< u_L< u_R<0$ 
or  $ u_L< u_b< u_R<0$ 
or  $ u_L< u_b<0< u_R$ 
or  $ u_L<0< u_b< u_R, \,\,\, u_b + u_L<0.$ In this case
$(u_L, \rho_L)$ is connected to $(u_R, \rho_R)$ by a rarefaction and the solution  is, 
\begin{equation}
\begin{aligned}
( u (x, t), \rho(x, t))= \begin{cases}
(u_L, \rho_L),  &\text{ if }   0< x< u_L t + x_0 ,\,\,\,  \forall t, \\
( \frac{x - x_0}{t}, 0),  &\text{ if }   u_L t+ x_0 < x<  u_R t+ x_0, \,\,\,  \forall t, \\
(u_R, \rho_R),  &\text{ if }   x> u_R t + x_0,\,\,\,   \forall t.
\end{cases}
\end{aligned}
\label{e3.23}
\end{equation}

Subcase 6:  If $ u_b< u_R< u_L<0.$ Here
 $(u_L, \rho_L)$ is connected to $(u_R, \rho_R)$ by a shock and the solution  is, 
\begin{equation*}
( u (x, t), \rho(x, t))
\end{equation*}
\begin{equation}
\begin{aligned}
= \begin{cases}
(u_L, \rho_L),  &\text{ if }   0< x< (\frac{ u_R + u_L}{2}) t + x_0, \\
( \frac{1}{2} (u_L+u_R), \frac{1}{2} (u_R- u_L)(\rho_R+\rho_R) t \delta_{x= (\frac{ u_L + u_R }{2}) t}),   &\text{ if }    x=( \frac{ u_R + u_L}{2}) t + x_0, \\
(u_R, \rho_R),  &\text{ if }    x>  (\frac{ u_R + u_L}{2}) t + x_0.
\end{cases}
\end{aligned}
\label{e3.24}
\end{equation}

Subcase 7: If $u_b <0< u_R < u_L$ or $ u_b <u_R<0< u_L.$ (Refer Figure 8.) In this case a rarefaction originating at $(0,0)$ interact with a shock originating at $(x_0,0)$ initially at $(x_1, t_1)$, which is the point of intersection of
$x =u_L t$ and $x =  \frac{ u_R + u_L}{2} t +x_0$.  The shock curve for $t > t_1$ is given by 
\begin{equation*}
\frac{dx}{dt}= \frac{1}{2} ( \frac{x}{t} + u_R), \,\,\, x(t_1)= x_1,
\end{equation*}
that gives the curve $x= \gamma_2(t)$ i.e.
$x= \frac{u_R}{2} t + c {t}^{ \frac{1}{2}}, \,\,\, c =  x_1 {t_1}^{-\frac{1}{2}} -  \frac{u_R}{2} {t_1}^{ \frac{1}{2}}$. The solution is
 
\begin{equation*}
( u (x, t), \rho(x, t))
\end{equation*}
\begin{equation}
\begin{aligned}
= \begin{cases}
( \frac{x}{t}, 0), &\text{ if }0< x<  u_L t,  \,\,\, t<t_1 \\
(u_L, \rho_L),  &\text{ if } u_L t< x<(\frac{ u_R + u_L}{2}) t +x_0, \,\,\,    t<t_1 \\
( \frac{1}{2} (u_L+u_R), \frac{1}{2} (u_R- u_L)(\rho_R+\rho_L) t \delta_{x= (\frac{ u_L + u_R }{2}) t}),  &\text{ if }  x= (\frac{ u_R + u_L}{2}) t +x_0,   \,\,\, t<t_1 \\
(u_R, \rho_R),  &\text{ if }  x> (\frac{ u_R + u_L}{2}) t +x_0,   \,\,\, t<t_1 \\
( \frac{x}{t}, 0),  &\text{ if } 0< x< \gamma_2(t), \,\,\,  t>t_1 \\
e(t) {\delta}_{x=\gamma_2(t)},  &\text{ if } t > t_1, \,\,\, x= \gamma_2(t), \\
(u_R, \rho_R),  &\text{ if }  x>  \gamma_2(t),  \,\,\,  t>t_1.
\end{cases}
\end{aligned}
\label{e3.25}
\end{equation} 
where shock strength has the form $e(t)=a t +b t^{1/2}+c$ where $a,b, c$ are constants and obtained by solving \eqref{e3.3} with initial condition $e(t_1)= \frac{1}{2} (u_R- u_L)(\rho_L+\rho_R) t_1$.

Subcase 8: If $0< u_b < u_R < u_L.$ (Refer Figure 8.) In this case $(u_b, \rho_b)$ is connected to $(u_L, \rho_L)$ by a rarefaction and $(u_L, \rho_L)$ is connected to $(u_R, \rho_R)$ by a shock and the interact initially at $(x_1,t_1)$ determined by the  point of intersection of  $x =u_L t$ and $x =  \frac{ u_R + u_L}{2} t +x_0$. After interaction the curve for $t > t_1$ is given by 
\begin{equation*}
\frac{dx}{dt}= \frac{1}{2} ( \frac{x}{t} + u_R), \,\,\, x(t_1)= x_1,
\end{equation*}
that gives the curve $x= \gamma_3(t)$ i.e.
$x= \frac{u_R}{2} t + c {t}^{ \frac{1}{2}}, \,\,\, c =  x_1 {t_1}^{-\frac{1}{2}} -  \frac{u_R}{2} {t_1}^{ \frac{1}{2}}$.
  The solution in this case is, 
\begin{equation*}
( u (x, t), \rho(x, t))
\end{equation*}
\begin{equation}
\begin{aligned}
= \begin{cases}
(u_b, \rho_b),  &\text{ if } 0< x<u_b t, \,\,\,  \forall t, \\
( \frac{x}{t}, 0), &\text{ if }u_b t< x<  u_L t,  \,\,\, t<t_1, \\
(u_L, \rho_L),  &\text{ if } u_L t< x<(\frac{ u_R + u_L}{2}) t +x_0, \\
 & \hspace{0.5cm} \text{and }     t<t_1, \\
( \frac{1}{2} (u_L+u_R), \frac{1}{2} (u_R- u_L)(\rho_R+\rho_L) t \delta_{x= (\frac{ u_L + u_R }{2}) t}),  &\text{ if }  x= (\frac{ u_R + u_L}{2}) t +x_0,   \,\,\, t<t_1 \\
(u_R, \rho_R),  &\text{ if }  x> (\frac{ u_R + u_L}{2}) t +x_0,   \,\,\, t<t_1, \\
( \frac{x}{t}, 0),  &\text{ if } u_b t< x< \gamma_3(t), \,\,\,  t>t_1, \\
e(t) {\delta}_{x=\gamma_2(t)},  &\text{ if } t > t_1, \,\,\, x= \gamma_2(t), \\
(u_R, \rho_R),  &\text{ if }  x>  \gamma_3(t),  \,\,\,  t>t_1.
\end{cases}
\end{aligned}
\label{e3.26}
\end{equation}
where shock strength has the form $e(t)=a t +b t^{1/2}+c$ where $a,b, c$ are constants and obtained by solving \eqref{e3.3} with initial condition $e(t_1)= \frac{1}{2} (u_R- u_L)(\rho_L+\rho_R) t_1$.

Subcase 9: If $0< u_L < u_b < u_R$ 
or $ u_L<0 < u_b < u_R, \,\,\, u_b + u_L >0$. (Refer Figure 9.) In this case $(u_b, \rho_b)$ is connected to $(u_L, \rho_L)$ by a shock and $(u_L, \rho_L)$ is connected to $(u_R, \rho_R)$ by a rarefaction and they interact. The initial point of interaction is $(x_1, t_1)$, is the point of intersection of the curves
 $x =u_L t+ x_0$ and $x =  \frac{ u_b + u_L}{2} t$ intersect at say . The shock curve for $t > t_1$ is given by 
\begin{equation*}
\frac{dx}{dt}= \frac{1}{2} ( \frac{x}{t} + u_b), \,\,\, x(t_1)= x_1,
\end{equation*}
that gives the curve $x= \gamma_4(t)$ i.e.
$x= \frac{u_b}{2} t + c {t}^{ \frac{1}{2}}, \,\,\, c =  x_1 {t_1}^{-\frac{1}{2}} -  \frac{u_b}{2} {t_1}^{ \frac{1}{2}}$.
 The solution in this case is, 
\begin{equation*}
( u (x, t), \rho(x, t))
\end{equation*}
\begin{equation}
\begin{aligned}
= \begin{cases}
(u_b, \rho_b),  &\text{ if } 0< x< (\frac{ u_b + u_L}{2}) t, \,\,\,  t< t_1, \\
( \frac{1}{2} (u_L+u_b), \frac{1}{2} (u_L- u_b)(\rho_b+\rho_L) t \delta_{x= (\frac{ u_L + u_b }{2}) t}),  &\text{ if }  x= (\frac{ u_b + u_L}{2}) t, \,\,\,  t< t_1, \\
(u_L, \rho_L), &\text{ if } (\frac{ u_b + u_L}{2}) t< x<  u_L t +x_0, \,\,\,    t<t_1, \\
( \frac{x - x_0}{t}, 0),  &\text{ if } u_L t + x_0 < x< u_R t + x_0, \,\,\,    t<t_1, \\
(u_R, \rho_R),  &\text{ if }  x>  u_R t + x_0,   \,\,\, \forall t, \\
(u_b, \rho_b),  &\text{ if } 0< x< \gamma_4(t), \,\,\,  t> t_1, \\
e(t) {\delta}_{x=\gamma_4(t)},  &\text{ if } t > t_1, \,\,\, x= \gamma_4(t), \\
( \frac{x - x_0}{t}, 0),  &\text{ if }  \gamma_4(t)< x< u_R t +x_0, \,\,\,   t> t_1.
\end{cases}
\end{aligned}
\label{e3.27}
\end{equation}
where shock strength has the form $e(t)=a t +b t^{1/2}+c$ where $a,b, c$ are constants and obtained by solving \eqref{e3.3} with initial condition $e(t_1)= \frac{1}{2} (u_L- u_b)(\rho_L+\rho_b) t_1$. 

Subcase 10: If $0< u_R < u_L < u_b$ or  $ u_R<0 < u_L < u_b$ 
or  $ u_R < u_L<0 < u_b, u_L + u_b>0.$ (Refer Figure 10.) In this case the shock curves originating at $(0,0)$ and $(x_0,0)$
namely the curves $x = (\frac{ u_b + u_L}{2}) t$ and $x =  (\frac{ u_R + u_L}{2}) t+ x_0$ intersect at say $(x_1, t_1)$.  In $t<t_1, \,\,\, (u_b, \rho_b)$ is connected to $(u_L, \rho_L)$ by a shock and $(u_L, \rho_L)$ is connected to $(u_R, \rho_R)$ by a shock. In $t>t_1,\,\,\, (u_b, \rho_b)$ is connected to $(u_R, \rho_R)$ by a shock. If $u_b+ u_R>0,$ then the curve   $x =  (\frac{ u_b + u_R}{2}) t+ x_1$ does not intersect the $t-$ axis at $t> t_1$. Whereas if $u_b+ u_R<0,$ then the curve   $x = ( \frac{ u_b + u_R}{2}) t+ x_1$ intersects the $t-$ axis at  $(0, \frac{-2 x_1}{ u_b + u_R})$ at $t> t_1$. The solution in both these cases is, 
\begin{equation*}
( u (x, t), \rho(x, t))
\end{equation*}
\begin{equation}
\begin{aligned}
= \begin{cases}
(u_b, \rho_b),  &\text{ if } 0< x< (\frac{ u_b + u_L}{2}) t, \,\,\,  t< t_1, \\
( \frac{1}{2} (u_L+u_b), \frac{1}{2} (u_L- u_b)(\rho_b+\rho_L) t \delta_{x= (\frac{ u_L + u_b }{2}) t}),  &\text{ if }  x= (\frac{ u_b + u_L}{2}) t, \,\,\,  t< t_1, \\
(u_L, \rho_L), &\text{ if } (\frac{ u_b + u_L}{2}) t< x< (\frac{ u_R + u_L}{2}) t +x_0, \,\,\, t<t_1, \\
(u_R, \rho_R),  &\text{ if }  x> (\frac{ u_R + u_L}{2}) t +x_0, \,\,\,  t<t_1, \\
(u_b, \rho_b),  &\text{ if } 0< x< (\frac{ u_R + u_b}{2}) t +x_1,  \,\,\,  t> t_1, \\
( \frac{1}{2} (u_R+u_b), \frac{1}{2} (u_R- u_b)(\rho_b+\rho_R) t \delta_{x= (\frac{ u_R + u_b }{2}) t}),  &\text{ if }  x= (\frac{ u_b + u_R}{2}) t + x_1, \,\,\,  t > t_1, \\
(u_R, \rho_R),  &\text{ if }  x>  (\frac{ u_R + u_b}{2}) t +x_1, \,\,\, t> t_1.
\end{cases}
\end{aligned}
\label{e3.28}
\end{equation}

Subcase 11: If $u_R < u_L < u_b<0$ 
or  $ u_R < u_L<0 < u_b, u_L + u_b<0$ 
or  $ u_R < u_b< u_L<0.$  In this case
  $(u_L, \rho_L)$ is connected to $(u_R, \rho_R)$ by a shock. This shall meet at $(0, t_1)$ on $x=0.$ The solution in this case is, 
\begin{equation*}
( u (x, t), \rho(x, t))
\end{equation*}
\begin{equation}
\begin{aligned}
= \begin{cases}
(u_L, \rho_L),  &\text{ if } 0< x< (\frac{ u_R + u_L}{2}) t+ x_0,  t< t_1, \\
( \frac{1}{2} (u_L+u_R), \frac{1}{2} (u_R- u_L)(\rho_L+\rho_R) t \delta_{x= (\frac{ u_L + u_R }{2}) t}),  &\text{ if }  x= (\frac{ u_R + u_L}{2}) t+ x_0, \,\,\,   t < t_1, \\
(u_R, \rho_R), &\text{ if } x>  (\frac{ u_R + u_L}{2}) t+ x_0,  \,\,\, x>0.
\end{cases}
\end{aligned}
\label{e3.29}
\end{equation} 

\begin{figure}[!hbtp]
\includegraphics[width=14cm,height=7cm]{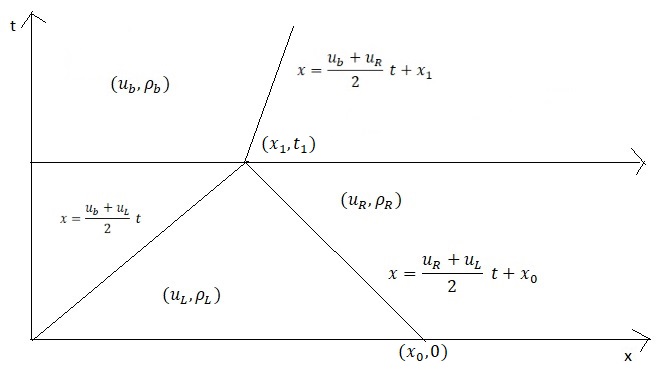}
\caption{ Interaction of shock centered at $(0,0)$ and shock centered at $(x_0,0)$}
\end{figure}

Subcase 12: If $0< u_R < u_b < u_L$ 
or $ u_R<0 < u_b < u_L.$ (Refer Figure 8.) Here a rarefaction centered at $(0,0)$ connecting $(u_b,\rho_b)$ to $(u_L,\rho_L)$ interact with a shock starting at $(x_0,0)$  connecting $(u_L,\rho_L)$ to $(u_R,\rho_R)$ with initial interaction point $(x_1,t_1)$ which is the point of intersection of the characteristic curve
 $x =u_L t$ and shock curve $x = ( \frac{ u_R + u_L}{2}) t+x_0$.  The shock curve for $t > t_1$ is given by 
\begin{equation*}
\frac{dx}{dt}= \frac{1}{2} ( \frac{x}{t} + u_R), \,\,\, x(t_1)= x_1,
\end{equation*}
that gives the curve $x= \gamma_5(t)$ i.e.
$x= \frac{u_R}{2} t + c {t}^{ \frac{1}{2}}, \,\,\, c =  x_1 {t_1}^{-\frac{1}{2}} -  \frac{u_R}{2} {t_1}^{ \frac{1}{2}}$.
 The solution in this case is, 
\begin{equation*}
( u (x, t), \rho(x, t))
\end{equation*}
\begin{equation}
\begin{aligned}
= \begin{cases}
(u_b, \rho_b),  &\text{ if } 0< x<u_b t, \,\,\,  t< t_1, \\
( \frac{x}{t}, 0),  &\text{ if } u_b t < x< u_L t, \,\,\,  t<t_1, \\
(u_L, \rho_L), &\text{ if } u_L t< x<  (\frac{ u_R + u_L}{2}) t+x_0,   t<t_1, \\
( \frac{1}{2} (u_L+u_R), \frac{1}{2} (u_R- u_L)(\rho_L+\rho_R) t \delta_{x= (\frac{ u_L + u_R }{2}) t}),  &\text{ if }  x= (\frac{ u_R + u_L}{2}) t+ x_0, \,\,\,     t<t_1, \\
(u_R, \rho_R),  &\text{ if }  x> (\frac{ u_R + u_L}{2}) t+x_0, \,\,\,  t<t_1, \\
( \frac{x}{t}, 0),  &\text{ if } u_b t< x< \gamma_5(t), \,\,\,  t> t_1, \\
e(t) {\delta}_{x=\gamma_5(t)},  &\text{ if } t > t_1, \,\,\, x= \gamma_5(t), \\
(u_R, \rho_R),  &\text{ if } x> \gamma_5(t), \,\,\, t> t_1.
\end{cases}
\end{aligned}
\label{e3.30}
\end{equation} 
where shock strength has the form $e(t)=a t +b t^{1/2}+c$ where $a,b, c$ are constants and obtained by solving \eqref{e3.3} with initial condition $e(t_1)= \frac{1}{2} (u_R- u_L)(\rho_L+\rho_R) t_1$. 

Subcase 13: If $u_R < u_b<0 < u_L.$ (Refer Figure 8 with $u_b=0$.)Here again the rarefaction centered at $(0,0)$ interact with a shock starting at $(x_0,0)$ with initial interaction point $(x_1,t_1)$ which is the point of intersection of the characteristic curve
 $x =u_L t$ and shock curve $x = ( \frac{ u_R + u_L}{2}) t+x_0$.  The shock curve for $t > t_1$ is given by  
\begin{equation*}
\frac{dx}{dt}= \frac{1}{2} ( \frac{x}{t} + u_R), \,\,\, x(t_1)= x_1,
\end{equation*}
that gives the curve $x= \gamma_6(t)$ i.e.
$x= \frac{u_R}{2} t + c {t}^{ \frac{1}{2}}, \,\,\, c =  x_1 {t_1}^{-\frac{1}{2}} -  \frac{u_R}{2} {t_1}^{ \frac{1}{2}}$.
The solution  in this case is, 
\begin{equation*}
( u (x, t), \rho(x, t))
\end{equation*}
\begin{equation}
\begin{aligned}
= \begin{cases}
( \frac{x}{t}, 0),  &\text{ if } 0< x< u_L t, \,\,\,  t<t_1, \\
(u_L, \rho_L), &\text{ if } u_L t< x<  (\frac{ u_R + u_L}{2}) t+x_0,  \,\,\, t<t_1, \\
( \frac{1}{2} (u_L+u_R), \frac{1}{2} (u_R- u_L)(\rho_L+\rho_R) t \delta_{x= (\frac{ u_L + u_R }{2}) t}),  &\text{ if }  x= (\frac{ u_R + u_L}{2}) t+ x_0, \,\,\,     t<t_1, \\
(u_R, \rho_R),  &\text{ if }  x> (\frac{ u_R + u_L}{2}) t+x_0, \,\,\,  t<t_1, \\
( \frac{x}{t}, 0),  &\text{ if } 0< x< \gamma_6(t), \,\,\,  t> t_1, \\
e(t) {\delta}_{x=\gamma_6(t)},  &\text{ if } t > t_1, \,\,\, x= \gamma_6(t), \\
(u_R, \rho_R),  &\text{ if } x> \gamma_6(t), \,\,\, t> t_1.
\end{cases}
\end{aligned}
\label{e3.31}
\end{equation}
where shock strength has the form $e(t)=a t +b t^{1/2}+c$ where $a,b, c$ are constants and obtained by solving \eqref{e3.3} with initial condition $e(t_1)= \frac{1}{2} (u_R- u_L)(\rho_L+\rho_R) t_1$.

\section{Conclusion}
 In previous studies on \eqref{1.1}, different methods have been used to analyse the  the initial value problem to \eqref{1.1} for the inviscid system. The solution was constructed using the Lax formula in \cite{le1} and existence of several solutions to the  Riemann problem was proved. The modified adhesion approximation  \eqref{1.3} was used in \cite{j1} to construct the solution  that gave the accurate physical solution coinciding with the solution constructed in \cite{le1}.  In this paper we construct weak asymptotic solution of initial boundary value problem and analysed interaction of  two waves,  one boundary Riemann solution originating at $(0,0)$ and another Riemann solution with initial discontinuity at $(x_0,0)$. We constructed explicit formula for the solution after interactions. This interaction problem is important when we study general initial boundary value problem, when different approximations such as Glimm scheme or Godunov scheme are used for actual computation of the solution.

\section*{Acknowledgments}

Funding: No funding involved. \\

Conflicts of interest/Competing interests: None. \\

Availability of data and material: Not applicable. \\

Code availability: Not applicable.

\section*{ORCID}
Kayyunnapara Divya Joseph https://orcid.org/0000-0002-4126-7882

\end{document}